\documentclass[a4paper,twoside,11pt,reqno]{amsart}
\usepackage[T1]{fontenc}
\usepackage{courier}
\usepackage[scaled=0.92]{helvet}

\usepackage{amscd,amsmath,amssymb,amsthm,amsfonts,epsfig,graphics}

\usepackage[english]{babel}
\usepackage{amsthm,amsfonts,amssymb,amsmath}
\usepackage{hyperref}
\hypersetup{
colorlinks=true,
citecolor=red,
linkbordercolor={1 1 1},
citebordercolor={1 1 1},
}

\usepackage{mathrsfs}

\usepackage{xcolor} 
\theoremstyle{plain}
\newtheorem{satz}{Theorem}[section]
\newtheorem{prop}[satz]{Proposition}

\newtheorem{lem}[satz]{Lemma}

\theoremstyle{definition}

\newtheorem{rem}[satz]{Remark}

\newtheorem{hyp}[satz]{Hypothesis}

\newcommand{\la}{\langle}
\newcommand{\re}{\rangle}
\newcommand{\mx}{\mbox}
\newcommand{\rw}{\rightarrow}
\newcommand{\de}{\displaystyle}
\newcommand{\ml}{\mathcal}

\newcommand{\pl}{\partial}

\newcommand{\x}{\times}

\newcommand{\beq}[1]{\begin{equation} \label{#1}}
\newcommand{\eeq}{\end{equation}}
\newcommand{\beqar}{\[ \begin{array}{rcl}}
\newcommand{\eeqar}{\end{array} \]}

\newcommand{\tild}{\tilde{d}}

\providecommand{\ep}{\varepsilon}
\providecommand{\ph}{\varphi}

\providecommand{\RR}{\mathbb{R}}
\providecommand{\CC}{\mathbb{C}}

\providecommand{\ZZ}{\mathbb{Z}}
\providecommand{\NN}{\mathbb{N}}
\providecommand{\TT}{\mathbb{T}}

\newcommand{\snorm}[2]{\left| #1\right|_{#2}}
\newcommand{\norm}[2]{\left \lVert#1 \right\rVert_{#2}}

\newcommand{\Heq}[2]{\overset{\left(#1\right)}{\underset{}{#2}}}

\usepackage[utf8]{inputenc}
\hypersetup{pdfauthor={Fortunati, Wiggins},%
            pdftitle={Aperiodic Kolmogorov Theorem},%
            pdfsubject={Hamiltonian Mechanics},%
            pdfkeywords={KAM Theorem, Aperiodic time dependence}%
}

\DeclareMathOperator{\id}{Id}

\makeatletter
\g@addto@macro{\endabstract}{\@setabstract}
\newcommand{\authorfootnotes}{\renewcommand\thefootnote{\@fnsymbol\c@footnote}}%
\makeatother

\begin{document}
%On the persistence of invariant tori under slowly decaying aperiodic perturbation
\title[Persistence of Diophantine flows...]{Persistence of Diophantine flows for quadratic nearly-integrable Hamiltonians under slowly decaying aperiodic time dependence}
\author{Alessandro Fortunati}
\thanks{This research was supported by ONR Grant No.~N00014-01-1-0769 and MINECO: ICMAT Severo Ochoa project SEV-2011-0087.}
\address{School of Mathematics, University of Bristol, Bristol BS8 1TW, United Kingdom}
\email{alessandro.fortunati@bristol.ac.uk}

\author{Stephen Wiggins}
\email{s.wiggins@bristol.ac.uk}

\maketitle

\begin{abstract}
The aim of this paper is to prove a Kolmogorov-type result for a nearly-integrable Hamiltonian, quadratic in the actions, with an aperiodic time dependence. The existence of a torus with a prefixed Diophantine frequency is shown in the forced system, provided that the perturbation is real-analytic and (exponentially) decaying with time. The advantage consists of the possibility to choose an arbitrarily small decaying coefficient, consistently with the perturbation size.
\\ 
The proof, based on the Lie series formalism, is a generalization of a work by A. Giorgilli.
\end{abstract}
\section{Introduction} 
The celebrated Kolmogorov Theorem, stated in \cite{kolm} with a guideline for the proof, has been for years a fruitful
source of ideas, culminating in the collection of tools and techniques nowadays known as KAM theory. As undisputed members of the acronym, Arnold \cite{arn1} and Moser \cite{mos62}, \cite{moser67} proposed complete proofs of
Kolmogorov's result. The two approaches exhibited some technical differences, but were both based on the concepts of
\emph{super-convergent method} and \emph{implicit function theorem} over the complexified phase space (see e.g.
\cite{chierchia09} for a detailed exposition). The  applicability of these tools to certain infinite dimensional
problems were investigated in \cite{moser1966}, giving rise to the modern theory of Nash-Moser arguments (see
\cite{zehnder76} and \cite{berbolpro} for an advanced setting).\\
The proof based on the Lie formalism proposed in \cite{nuovocimento} then continued in \cite{giorgloca97}, \cite{giorgmorbi97} and \cite{giorgloca99},  makes use of the well known class of canonical change in \emph{explicit form}. This has the remarkable advantage to avoid the inversion and the difficulties related to implicit function arguments. Furthermore, this feature has been widely and profitably used for the computer implementation of  normalization algorithms.\\ 
In a substantially different direction, the approach developed in \cite{chierchiafalco94}, \cite{chierchiafalco96} and by the Gallavotti's school \cite{gallavottitwistlesskam}, \cite{gallagent95}, \cite{gentmastro95} and subsequent papers, is based on \emph{renormalization group} tools and  \emph{diagrammatic} analysis of the Lindstedt's series convergence due to cancellation phenomena. The analysis is an extensive improvement of the pioneering challenge of the small divisors problem faced in \cite{eliasson88}.\\
The historical legacy between the Kolmogorov Theorem and problems arising from  Celestial Mechanics, has led to a development in the treatment of quasi-periodic perturbations of integrable Hamiltonians, mainly in the presence of weaker regularity hypothesis.\\  
Our aim is to proceed in a slightly different direction, investigating the possibility of obtaining the conservation of (strongly) non-resonant tori in the case of an analytic perturbation (quadratic in the actions), but with an \emph{aperiodic} time dependence. For this purpose we shall follow the  exposition \cite{gior}, a revisited essay of the techniques used in \cite{nuovocimento}. The case of a quadratic Hamiltonian, has been chosen for simplicity of discussion. On the other hand, this choice allows substantial simplification of the ``known'' technical part, emphasizing the differences introduced by the non-quasi-periodic time dependence. As we shall discuss, the exponential rate of the perturbation decay, say $\exp(-at)$, is a simplified choice as well. \\
The philosophy behind the present analysis is very close to the Nekhoroshev stability result for aperiodically perturbed system of \cite{forwig}, but some substantial differences arise. Mainly, the Nekhoroshev normal form can be constructed by modifying the original normalization scheme, with the sole hypothesis that the perturbation depends $\mu-$slowly
on time. Hence the technical part consists in giving an estimate of the extra-terms arising from the aperiodic dependence. The key point is that, as it is clear from the normal form statement (see \cite[Thm 2]{forwig}), this is possible only because the number $r$ of normalization steps is \emph{finite} and the threshold for $\mu$ is actually a function of $r$.\\
The same phenomenon,
even in the presence of a different normalization scheme, can be found if the Kolmogorov construction is extended
\emph{tout-court} to the case of aperiodic perturbations, and the slow dependence hypothesis would inevitably degenerate to a trivial case
 i.e. $\mu = 0$.\\
The above described difficulty, has required the modification of the transformation suggested by Kolmogorov in a way to annihilate certain time dependent terms arising in the normalization algorithm. The standard homological equation is modified, in this way, into a linear PDE involving time. The
apparently ``cheating'' hypothesis of time decaying perturbation (asymptotically the problem is trivial)
turns out to be a technical ingredient in order to ensure the resolvability of this equation at each step
of the normal form construction. Nevertheless, as a feature behind the \emph{slow decay}, the whole argument does not impose lower bounds on $a$. Consistently, the slower the decay, the smaller the perturbation size.\\
The self-contained exposition is closely carried along the lines of \cite{gior}. The same notational setting is used for a more efficient comparison.

\section{Preliminaries and statement of the result}
Let us consider the following Hamiltonian
\beq{eq:hamtempo}
\ml{H}(Q,P,t)=\frac{1}{2}\la \Gamma P,P \re + 
\ep f(Q,P,t) \mx{,}
\eeq
where $\Gamma$ is a $n \times n$ real symmetric matrix, $(Q,P) \in \TT^n \times \RR^n$ is a set of action-angle variables, $t \in \RR^+$ is an additional variable (time) and $\ep>0$ is a small parameter.  The perturbing function $f$ is assumed to be quadratic in $P$. \\
The Kolmogorov approach to (\ref{eq:hamtempo}) begins by considering a given $\hat{P} \in \RR^n$ then expanding the first term of $\ml{H}$ around it. The canonical change (translation) $(q,p):=(Q,P-\hat{P})$, and the definition of $\eta \in \RR$ as the momentum conjugate to $\xi:=t$, yields (up to a constant) the following autonomous Hamiltonian
\beq{eq:ham}  
H(q,p,\xi,\eta):=
\la \omega, p \re + \frac{1}{2}\la \Gamma p,p \re  + \eta + \ep f(q,p,\xi) \mx{,}
\eeq
where $\omega:=\Gamma \hat{P}$.\\
In order to use the standard tools concerning analytic functions, we consider a complex extension of the ambient space. More precisely, define $\ml{D}:=\Delta_{\rho} \times \TT_{2 \sigma}^n \x \ml{S}_{\rho} \x \ml{R}_{\zeta}$ where
$$
\begin{array}{rclrcl}
\Delta_{\rho}&:=&\{p \in \CC^n:|p|<\rho\},
& \qquad 
\TT_{2 \sigma}^n&:=&\{q \in \CC^n: |\Im q| < 2 \sigma\},\\
\ml{S}_{\rho}&:=&\{\eta \in \CC: |\Im \eta| <\rho\}, & \qquad 
\ml{R}_{\zeta}&:=&\{\xi =:x+iy \in \CC:|x| <\zeta\ ; \,
y > -\zeta\} \mx{,}
\end{array}
$$
and $\rho,\sigma,\zeta \in (0,1)$. Similarly to \cite{gior}, we consider the usual \emph{supremum norm}
$$
\snorm{g}{[\rho,\sigma;\zeta]}:=\sup_{(p,q) \in \ml{D}} |g(q,p,\xi)| \mx{,}
$$
and the \emph{Fourier norm}, defined for all $\nu \in(0,1/2]$,
\beq{eq:fouriernorm}
\norm{g}{[\rho,\sigma;\zeta]}:=\sum_{k \in \ZZ^n} \snorm{g_k(p,\xi)}{(\rho,\sigma)} 
e^{2|k|(1-\nu)\sigma} \mx{,}
\eeq
where $g_k(p,\xi)$ are the coefficient of the Fourier expansion $g=\sum_{k \in \ZZ^n} g_k(p,\xi) e^{i \la k , q \re}$. For all vector-valued functions $w:\ml{D} \rw \CC^n $ we shall set $\norm{w}{[\rho,\sigma;\zeta]}:=\sum_{l=1}^n \norm{w_l}{[\rho,\sigma;\zeta]}$.\\
System (\ref{eq:ham}) will be studied under the following 
\begin{hyp}\label{hyp}
\begin{itemize}
\item There exists $m\in (0,1)$ such that, for all $ v \in \CC^n$ 
\beq{eq:hypongamma}
 |\Gamma v| \leq m^{-1}|v| \mx{.}
\eeq
\item (Slow decay): The perturbation is an analytic function on $\ml{D}$ satisfying
\beq{eq:slowdecay}
\norm{f(q,p,\xi)}{[\rho,\sigma;\zeta]} \leq M_f e^{-a |\xi|} \mx{,}
\eeq
for some $M_f>0$ and $a \in (0,1)$. 
\end{itemize}
\end{hyp}
We specify that the assumption $a<1$ (which includes, of course, the ``interesting'' case of $a$ small) is not of technical nature, but it is often useful to obtain more compact estimates. As a difference with \cite{forwig}, hypothesis (\ref{eq:slowdecay}) is not of slow time dependence: in principle, the constant $M_f$ could be the bound of an arbitrary (analytic) function of $\ph$ and of $\xi$.\\
In this framework, the main result is stated as follows
\begin{satz}[Aperiodic Kolmogorov]\label{thm}
Consider Hamiltonian (\ref{eq:ham}) under the Hypothesis \ref{hyp} and suppose that $\hat{P}$ is such that $\omega$ is a $\gamma-\tau$ Diophantine vector\footnote{Namely, there exist $\gamma$ and $\tau>n-1$ such that $|\la \omega, k \re | \geq \gamma |k|^{-\tau}$, for all $k \in \ZZ^{n}\setminus \{0\}$, understood $|k|:=|k_1|+\ldots + |k_n|$.}.\\
Then, for all $a \in (0,1)$ there exists\footnote{See (\ref{eq:finalvalueep}) for an explicit estimate.} $\ep_a>0$ such that, for all $\ep \in (0,\ep_a]$,  
it is possible to find a canonical, $\ep-$close to the identity, analytic change of variables  $(q,p,\xi,\eta) = \ml{K}(q^{(\infty)},p^{(\infty)},\xi,\eta^{(\infty)})$, $\ml{K}:\ml{D}_* \rw \ml{D}$ with $\ml{D}_* \subset \ml{D}$, 
casting Hamiltonian (\ref{eq:ham}) into the \emph{Kolmogorov normal form} 
\beq{eq:kolnormal}
H_{\infty}(q^{(\infty)},p^{(\infty)},\xi,\eta^{(\infty)})=\la \omega, p^{(\infty)} \re + \eta^{(\infty)}+
\mathcal{Q}(q^{(\infty)},p^{(\infty)},\xi;\ep) \mx{,}
\eeq
with $\pl_p^{\alpha} \ml{Q}(\cdot,0,\cdot;\ep)=0$ for all $\alpha \in \NN^n$ such that $\alpha_i \leq 1$ ($\ml{Q}$ is a homogeneous polynomial of degree $2$ in $p$). 
\end{satz}
Hamiltonian (\ref{eq:kolnormal}) is defined up to a function of $\xi$ that is not relevant for the $(q,p)-$ flow we are interested in. The normal form (\ref{eq:kolnormal}) clearly implies the persistence of the (lower dimensional for (\ref{eq:ham}) i.e. maximal for (\ref{eq:hamtempo})) invariant torus with frequency $\omega$ under perturbations satisfying (\ref{eq:slowdecay}) and for sufficiently small $\ep$. \\
The rest of the paper is devoted to the proof of Theorem \ref{thm}. As usual, it has the structure of an iterative statement divided into a formal part (Lemma \ref{lem:iterativeformal}) and a quantitative part (Lemma \ref{lem:iterative}). In the first part we modify the Kolmogorov scheme in order to build a suitable normalization algorithm for the problem at hand. The  homological equation on $\TT_{2 \sigma}^n \times \ml{R}_{\zeta}$ arising in this case requires a substantially different treatment of the bounds on the small divisors as described in Prop. \ref{prop:small}.\\
In the second, quantitative part, the well established tools of the Lie series theory (recalled in Sec. \ref{sec:technical}), are used to control the size of the unwanted terms during the normalization process, proving that the constructed  Kolmogorov transformation has the feature to make them smaller and smaller.\\
The final part consists in showing that the described iterative argument can be iterated infinitely many times, and the contribution of the unwanted terms completely removed: once more, the choice of a particular torus $P=\hat{P}$ suggested by Kolmogorov, is required for the convergence of this particular scheme.

\section{The formal perturbative setting}
Following \cite{gior} we construct a perturbative scheme in which the $j-$th step is based on the canonical transformation 
\[
\ml{K}_j:=\exp(\ml{L}_{\chi^{(j)}}) \circ \exp(\ml{L}_{\phi^{(j)}}) \mx{,}
\]
where the \emph{Lie series operator} is formally defined by
\[
\exp(\ml{L}_{G}):=\id+\sum_{s \geq 1} \frac{1}{s!}\ml{L}_G^s \mx{,}
\]
and $\ml{L}_G \cdot :=\{G,\cdot\}=(\pl_q G \pl_p + \pl_{\xi}G \pl_{\eta}-\pl_p G \pl_q-\pl_{\eta}G \pl_{\xi})\cdot $ is the \emph{Lie derivative}. The \emph{generating functions} will be chosen of the form $\phi^{(j)}=\phi^{(j)}(q,\xi)$ and $\chi^{(j)}=\chi^{(j)}(q,p,\xi):=\la Y^{(j)}(q,\xi),p \re $. The latter being the equivalent of the classical case.
\begin{lem}\label{lem:iterativeformal}
Suppose that for some $j \in \NN$, Hamiltonian (\ref{eq:ham}) can be written in the form
\beq{eq:hamricorsiva}
H_j=\la \omega,p \re + \eta + A^{(j)}(q,\xi) + \la B^{(j)}(q,\xi),p \re
 + \frac{1}{2} \la C^{(j)}(q,\xi)p,p \re \mx{,}
\eeq
with $C^{(j)}$ symmetric matrix. Then it is possible to determine $\phi^{(j)}$ and $Y^{(j)}$ such that
$H_{j+1}:=\ml{K}_j H_j$ has the structure (\ref{eq:hamricorsiva})  for suitable $A^{(j+1)},B^{(j+1)}$ and $C^{(j+1)}$ symmetric matrix as well.
\end{lem} 
The possibility to write the Hamiltonian (\ref{eq:ham}) in the form (\ref{eq:hamricorsiva}), and then to complete an iterative scheme, will be discussed in Sec. \ref{seq:inducbasis}.
\begin{rem} 
The variables change casting $H_j$ into $H_{j+1}$ follows directly from the Gr{\"o}bner \emph{exchange} Theorem\footnote{Namely, let for simplicity $H=H(q,p)$ and $\chi$ be a generating function, one has
\[
H(q,p)|_{(q,p)=\exp(\ml{L}_{\chi})(q',p')}=[\exp(\ml{L}_{\chi})H(q,p)]_{(q,p)=(q',p')} \mx{,}
\]
understood $\exp(\ml{L}_{\chi})(q',p')=(\exp(\ml{L}_{\chi})q',\exp(\ml{L}_{\chi}) p')$.
} and reads as
\beq{eq:change}
(q^{(j)},p^{(j)},\eta^{(j)},\xi^{(j)})=\ml{K}_j
(q^{(j+1)},p^{(j+1)},\eta^{(j+1)},\xi^{(j+1)}) \mx{.}
\eeq
As a basic feature of this method, the variables superscript is not relevant in order to deal with the Hamiltonian transformation, and it will be omitted throughout the proof.
\end{rem}
The perturbative feature of this result is not transparent until a quantitative control of the action of $\ml{K}_j$ is established. Indeed, the subsequent step is to show that the ``size'' (in a sense that will be made precise later) of the terms $A^{(j)},B^{(j)}$ is infinitesimal as $j$ tends to infinity, obtaining in this way the desired \emph{Kolmogorov normal form}.  
\proof 
It is convenient to discuss separately the action of the two transformations.
\subsubsection*{First transformation} Firstly we examine the action of $\exp(\ml{L}_{\phi^{(j)}})$ on $H_j$. 
A key feature of $\ml{L}_{\phi^{(j)}}$, is that the degree of polynomials in $p$ on which it acts are decreased by one order. This implies that $\exp(\ml{L}_{\phi^{(j)}})H_j$ turns out to be simply 
\beqar
 \exp(\ml{L}_{\phi^{(j)}})H_j & = & \de \la \omega , p \re + \pl_{\omega} \phi^{(j)}  +  \eta +  \pl_{\xi} \phi^{(j)}  + A^{(j)}+  \la B^{(j)}, p \re+ \la B^{(j)}, \pl_q \phi^{(j)} \re\\
 & + &  \de   
 \frac{1}{2} \la C^{(j)} p,p \re +  \la C^{(j)} \pl_q \phi^{(j)}, p \re + 
 \frac{1}{2}\la C^{(j)} \pl_q \phi^{(j)} ,\pl_q \phi^{(j)} \re \mx{,}
 \eeqar 
 where $\pl_{\omega} \cdot :=\la \omega , \pl_q \cdot \re$. Note that the symmetry of $C^{(j)}$ has been repeatedly used. 
 \begin{rem} The finite number of terms in the previous expression is clearly one of the main simplifications introduced by a $p-$quadratic Hamiltonian. By considering the remainder of degree $\geq 3$ in $p$, the Lie series operator would have produced an infinite number of terms.
 \end{rem}
 The first generating function $\phi^{(j)}(q,\xi)$ is determined as the solution of the following \emph{time dependent homological equation}
\beq{eq:firsthomological}
\pl_{\xi} \phi^{(j)} (q,\xi)+ \pl_{\omega} \phi^{(j)} (q,\xi) + A^{(j)} (q,\xi)=0\mx{.}
\eeq
This equation can be formally solved on the Fourier space, giving rise to an infinite set of decoupled ODEs, see Prop. \ref{prop:small} for more details. In spite of this difficulty, the presence of the term $\pl_{\xi} \phi^{(j)}$ allows the resolvability of the equation also for the $0-$th Fourier coefficient ($q-$average\footnote{We shall denote also with $\overline{f(q,\xi)}:=(2 \pi)^{-n}\int_{\TT^n} f(q,\xi) dq$ the $q-$average of $f$.}). This feature, not necessary in this case ($\overline{A(q,\xi)}$ could be removed from this equation and kept in the Hamiltonian without affecting the normal form) will play a key role in the determination of $Y^{(j)}(q,\xi)$. Now, defining
\begin{subequations}
\begin{align}
\hat{A}^{(j)}(q,\xi)& :=  \de
\la  B^{(j)}, \pl_q \phi^{(j)} \re +  \frac{1}{2} \la C^{(j)} \pl_q \phi^{(j)} ,\pl_q \phi^{(j)} \re \mx{,}\label{eq:ahat}\\
\hat{B}^{(j)}(q,\xi)& :=  \de B^{(j)}+ C^{(j)} \pl_q \phi^{(j)}   \mx{,}\label{eq:bhat}
\end{align}
\end{subequations}
we obtain
\beq{eq:hamhalstep}
 \hat{H}_j:=\exp(\ml{L}_{\psi_j})H_j=\la \omega , p \re + \eta + \hat{A}^{(j)}(q,\xi)+ 
 \la \hat{B}^{(j)}(q,\xi),p \re + \frac{1}{2} \la C^{(j)}(q,\xi)p,p \re \mx{.}
\eeq
\subsubsection*{Second transformation} Our aim is now to determine $Y^{(j)}(q,\xi)$. Explicitly we have
\beqar
\exp(\ml{L}_{\chi^{(j)}}) \hat{H}_j & = & \id \hat{H}_j + \de \ml{L}_{\chi^{(j)}} \la \omega , p \re + \ml{L}_{\chi^{(j)}} \eta+ \sum_{s \geq 2} \frac{1}{s!} \ml{L}_{\chi^{(j)}}^s \la \omega , p \re+ \sum_{s \geq 1} \frac{1}{s!} \ml{L}_{\chi^{(j)}}^s \hat{A}^{(j)} \\
 & + & \de \sum_{s \geq 1} \frac{1}{s!} \ml{L}_{\chi^{(j)}}^s \la \hat{B}^{(j)},p \re + \de \sum_{s \geq 1} \frac{1}{s!} \ml{L}_{\chi^{(j)}}^s \la C^{(j)} p,p \re
+   \sum_{s \geq 2} \frac{1}{s!} \ml{L}_{\chi^{(j)}}^s \eta \mx{.}
\eeqar
The function $\chi^{(j)}(q,\xi)$ is determined in such a way 
\beq{eq:secondhompre}
\ml{L}_{\chi^{(j)}} \eta+\ml{L}_{\chi^{(j)}} \la \omega,p \re + \la \hat{B}^{(j)}(q,\xi),p \re=0 \mx{.}
\eeq
Noting that 
\beqar
\de \sum_{s \geq 2} \frac{1}{s!} \ml{L}_{\chi^{(j)}}^s \la \omega , p \re
+ \sum_{s \geq 1} \frac{1}{s!} \ml{L}_{\chi^{(j)}}^s \la \hat{B}^{(j)},p \re 
&=& \de \sum_{s \geq 1} \frac{1}{(s+1)!} \ml{L}_{\chi^{(j)}}^s 
[\ml{L}_{\chi^{(j)}} \la \omega , p \re + (s+1) \la \hat{B}^{(j)},p \re ] \\
& \Heq{\ref{eq:secondhompre}}{=}  & \de \sum_{s \geq 1} \frac{s}{(s+1)!} \ml{L}_{\chi^{(j)}}^s \la \hat{B}^{(j)},p \re 
-   \sum_{s \geq 2} \frac{1}{s!} \ml{L}_{\chi^{(j)}}^s \eta  \mx{,}
\\  
\eeqar
the transformed Hamiltonian simplifies as follows
\[
\exp(\ml{L}_{\chi_j}) \hat{H}_j  =  \de \la \omega , p \re +\eta + \exp(\ml{L}_{\chi^{(j)}})\hat{A}^{(j)}+ \sum_{s \geq 1} \frac{s}{(s+1)!} \ml{L}_{\chi^{(j)}}^s \la \hat{B}^{(j)},p \re + \frac{1}{2}\exp(\ml{L}_{\chi^{(j)}}) \la C^{(j)} p,p \re\mx{.}
\]
It is sufficient to define
\begin{subequations}
\begin{align}
A^{(j+1)}(q,\xi)  &:= \de \exp(\ml{L}_{\chi^{(j)}})\hat{A}^{(j)} \mx{,} \label{eq:ajpuno}\\
\la B^{(j+1)}(q,\xi),p \re  & :=   \de 
 \sum_{s \geq 1} \frac{s}{(s+1)!} \ml{L}_{\chi^{(j)}}^s
 \la \hat{B}^{(j)}, p \re  \mx{,} \label{eq:bjpuno}\\
\la C^{(j+1)}(q,\xi) p,p \re & :=  \exp(\ml{L}_{\chi^{(j)}}) \la C^{(j)} p,p \re \mx{,}
\label{eq:cjpuno}
\end{align}
\end{subequations}
in order to obtain
\beq{eq:transformed}
H_{j+1}:=\exp(\ml{L}_{\chi^{(j)}}) \hat{H}_j= \de \la \omega , p \re + \eta + 
 A^{(j+1)}(q,\xi)+ \la B^{(j+1)}(q,\xi),p \re + \frac{1}{2} \la C^{(j+1)}(q,\xi)p,p \re \mx{,}
\eeq
which has the structure (\ref{eq:hamricorsiva}). The symmetry of $C^{(j+1)}$ follows from its definition.\\
It is immediate to check that (\ref{eq:secondhompre}) is equivalent to $\la (\pl_{\xi} Y^{(j)}+ \pl_{\omega} Y^{(j)}  + \hat{B}^{(j)}), p \re=0$, i.e., 
\beq{eq:secondhom}
\pl_{\xi} Y^{(j)}(q,\xi)+ \pl_{\omega} Y^{(j)}(q,\xi) + \hat{B}^{(j)}(q,\xi)=0 \mx{,}
\eeq
which has the same form of (\ref{eq:firsthomological}) if considered component-wise. The necessity to solve (\ref{eq:secondhom}) also for the $0-$th Fourier mode is now clear: any ``residual'' term would imply a frequency correction and the failure of the program.

\section{Technical tools}\label{sec:technical}
From this section on, we shall profitably use the complex analysis tools in order to show the convergence of the Kolmogorov scheme. Let us firstly recall a well known property of the analytic functions: if  $g=g(q,p,\xi)$ is analytic on $\ml{D}$, one has $|g_k| \leq \snorm{g}{[\rho,\sigma;\zeta]} e^{-2|k|\sigma}$ then, by (\ref{eq:fouriernorm}), $\norm{g}{[\rho,\sigma;\zeta]}< \infty$ for all $\nu>0$. Vice-versa, if $\norm{g}{[\rho,\sigma;\zeta]}< \infty$ for all $\nu>0$ (no matter how small), then the Fourier coefficients of $g$ decay as $e^{-2|k|\sigma}$, hence the corresponding series defines an analytic function\footnote{I.e. the finiteness of the Fourier norm characterizes analytic functions on $\ml{D}$, see e.g. \cite[Chap. $4$]{giorgilli02}. The choice of $\nu$ will be tacitly understood in the follow as sufficiently small in order to ensure that the function at hand is analytic in a domain that is as large as possible.} on $\ml{D}$. \\
As in \cite{gior} we collect some basic inequalities in the following
\begin{prop} Let $v(q,\xi)$ and $C(q,\xi)$ respectively a vector and a matrix defined on $\ml{D}$. Then the following property hold
\begin{itemize}
\item 
\beq{eq:inequno}
\norm{\la v(q,\xi),p \re}{[\rho,\sigma;\zeta]} \leq \rho \norm{v}{[\sigma;\zeta]} \mx{.}
\eeq
Vice-versa, if for some $\tilde{M}>0$ 
\beq{eq:inequnocontrary}
\norm{\la v(q,\xi),p \re}{[\rho,\sigma;\zeta]} \leq \tilde{M} \rho, \qquad 
\mx{then} \qquad \norm{v(q,\xi)}{[\sigma;\zeta]} \leq \tilde{M} \mx{.} 
\eeq
\item If, for some $\hat{M}>0$ 
\beq{eq:ineqdue}
\norm{\la C(q,\xi)p,p \re }{[\rho,\sigma;\zeta]} \leq \hat{M} \rho^2, \qquad 
\mx{then} \qquad \norm{C_{kl}(q,\xi)}{[\sigma;\zeta]} \leq \hat{M} \mx{.}
\eeq
\end{itemize}
\end{prop}
\proof It can be extended without difficulties to our case, by following the sketch proposed in \cite[Pag. 160]{gior} \endproof
It will be also useful to recall the bound below, valid in particular on $\ml{R}_{\zeta}$ 
\beq{eq:ineq}
e^{-a|x|} \leq e^{a \zeta} e^{-a|\xi|} \mx{.}
\eeq 

\subsection{Solution of the time dependent homological equation} 
Let us consider the following P.D.E.
\beq{eq:hompde}
\pl_{\xi} \ph + \pl_{\omega} \ph = \psi \mx{,}
\eeq 
where $\psi=\psi(q,\xi):\mathcal{D} \rw \CC$ is a given function. It is possible to state the following
\begin{prop}\label{prop:small} Let $\delta \in [0,1)$ and suppose that $\psi$ is analytic on $\TT_{2 (1 - \delta) \sigma}^n \times \ml{R}_{\zeta} $ and \emph{exponentially decaying} with $|\xi|$, i.e. 
\beq{eq:expdec}
\norm{\psi}{[(1-\delta)\sigma;\zeta]} \leq K e^{-a |\xi|} \mx{,}
\eeq
where $a$ has been defined in (\ref{eq:slowdecay}).\\
Then for all $d \in (0, 1-\delta)$ and for all $\zeta$ such that
\beq{eq:sceltazeta}
2 |\omega| \zeta \leq d \sigma \mx{,}
\eeq
the solution of (\ref{eq:hompde}) exists and satisfies
\begin{subequations}
\begin{align}
\norm{\ph}{[(1-\delta-d)\sigma;\zeta]} & \leq  
\frac{K S_1}{a(d \sigma)^{\tau}} e^{-a |\xi|} \mx{,} \label{eq:stimehomuno}
\\
\norm{\pl_{q_m} \ph }{[(1-\delta-d)\sigma;\zeta]} & 
\leq \frac{K S_2}{a (d \sigma)^{\tau+1}} 
e^{-a |\xi|},\qquad m=1,\ldots,n \mx{,}
\label{eq:stimehomdue}
\end{align}
\end{subequations}
\end{prop}
where $S_{1,2}>0$ are constants defined for all sufficiently small $\nu>0$.
\proof 
By expanding $\ph=\ph(q,\xi)$ we have that equation (\ref{eq:hompde}) in terms of Fourier coefficients reads as 
\[
i \lambda \ph_k(\xi) + \ph_k'(\xi)= \psi_k(\xi) \mx{,}
\]
with $\lambda:=\la \omega , k \re$. We firstly discuss the case $k \neq 0$, hence $\lambda \neq 0$ by assumption. The solution in this case is
\[
\ph_k(\xi)=e^{-i \lambda \xi} \left[ \ph_k(0) + \int_0^\xi \psi_k(s) e^{i \lambda s} ds \right] \mx{.}
\]
The integral is meant to be computed along an arbitrary path ($\ml{R}_{\zeta}$ is simply connected) joining the origin and $\xi \in \CC$. More precisely,  we shall choose
\beq{eq:integralasymp}
\int_0^\xi \psi_k(s) e^{i \lambda s} ds=\int_0^x \psi_k(x') e^{i \lambda x'} dx' 
+ i e^{i \lambda x} \int_0^y \psi_k(x+i y') e^{-\lambda y'} dy' \mx{.} 
\eeq
The complex number $\ph_k(0)$ denotes the value of the solution at the complex plane origin and it will be determined in such a way 
$\lim_{\Re(\xi) \rw \infty} \ph_k(\xi)=0$, i.e. taking into account the hypothesis (\ref{eq:expdec})  
\[
\ph_k(0)=-\int_0^{+\infty} \psi_k(x) e^{i \lambda x} dx 
\mx{.}
\]
As a consequence, the solution satisfies
\[
|\ph_k(\xi)| \leq e^{\lambda y} \left[ 
\int_0^y |\psi_k(x+i y')| e^{-\lambda y'} dy'
+ 
\int_x^{+\infty} |\psi_k(x')| dx' 
\right] \mx{.}
\]
By hypothesis (\ref{eq:expdec}) it follows that $|\psi_k(\xi)| \leq K e^{-[a|\xi|+2|k|(1-\delta)\sigma]} $, hence the integrals appearing in the previous formula can be bounded on the strip $\ml{R}_{\zeta}$ as follows 
\beqar
\de \int_0^y |\psi_k(x+i y')| e^{-\lambda y'} dy' & \leq & \de K e^{-[a |x| +2|k|(1-\delta)\sigma] } \int_0^y e^{|\lambda|y'} dy' \\ [8pt]
& \leq & \de |\lambda|^{-1}K  e^{-[a |x| +2|k|(1-\delta)\sigma-|\lambda|\zeta]} \mx{,}\\ [5pt]
\de \int_x^{\infty}  |\psi_k(x')| dx' & \leq & \de K e^{-2|k|(1-\delta)\sigma} \int_x^{\infty} e^{-a |x'|} dx'\\ [8pt]
& \leq & \de 2 K a^{-1} e^{a \zeta} e^{-[a |x| +2|k|(1-\delta)\sigma]} \mx{.}
\eeqar
The obtained estimates imply
\beq{eq:stimafinalephk}
|\ph_k(\xi)| \leq 
K e^{-[a x +2|k|(1-\delta)\sigma-2|\lambda|\zeta]} \left[ \frac{1}{|\lambda|}+\frac{2e^{a \zeta}}{a}\right] \leq
2 K \frac{(a \gamma+e^{a \zeta})}{a} |k|^{\tau} e^{-[a |x| +2|k|(1-\delta)\sigma-2|\lambda|\zeta]}  \mx{,}
\eeq
where we used the Diophantine condition. Now using inequalities $|\lambda| \leq |k||\omega|$,  
\[
|k|^{\tau} e^{-d |k| \sigma} \leq \left( \frac{ \tau}{e d \sigma}\right)^{\tau} \mx{,}
\]
and finally hypothesis (\ref{eq:sceltazeta}), one has
\beq{eq:phk}
|\ph_k(\xi)| \leq 2 K \frac{(a \gamma+e^{a \zeta})}{a} 
\left( \frac{ \tau}{e d \sigma}\right)^{\tau}
e^{-a|x|} e^{-2|k|(1-\delta-d)\sigma}  \mx{.} 
\eeq
Hence the series $\sum_{k \in \ZZ^n \setminus \{0\}} \ph_k(\xi)$ defines an analytic function on $\TT_{2 (1 - \delta- d) \sigma}^n \times \ml{R}_{\zeta} $.\\ The simpler case $k=0$, yielding the equation $\pl_{\xi} \ph_0(\xi)=\psi_0(\xi)$,  can be treated in similar way. More precisely, by determining $\ph_0(0)$ as in (\ref{eq:integralasymp}) and bounding the two resulting integrals of the path we get
\beq{eq:phz}
|\ph_0(\xi)| \leq \zeta K e^{-a|x|} + \frac{2Ke^{a \zeta}}{a} e^{-a|x|} \leq \frac{4Ke^{a \zeta}}{a} e^{-a|x|} \mx{.}
\eeq
Now recall definition (\ref{eq:fouriernorm}). By (\ref{eq:phk}) and (\ref{eq:phz}), the use of (\ref{eq:ineq}) (recalling $a,\zeta<1$), and finally by setting
\[
S_1:=4 e^2+2 (\gamma+e)(\tau/e)^{\tau} \sum_{k \in \ZZ^n  \setminus \{0\}} e^{-2 \nu|k|(1-\delta-d)\sigma} 
\] 
we get (\ref{eq:stimehomuno}). Note that, as long as $d+\delta<1$, the upper bound for $S_1$ is independent on $d,\delta$, being $\nu$ arbitrarily small.  As for as $\pl_{q_m}\ph$, directly from the Fourier expansion we find
$
\pl_{q_m} \ph(q,\xi)=i \sum_{k \in \ZZ^n\setminus\{0\}} k_m \ph_k(\xi)e^{i \la k , q \re}$. By using bound (\ref{eq:stimafinalephk}) (the average term is not relevant in such case) and proceeding in a similar way we get (\ref{eq:stimehomdue}), where $S_2:=[(\tau+1)/e]^{(\tau+1)}\sum_{k \in \ZZ^n  \setminus \{0\}} e^{-2 \nu|k|(1-\delta-d)\sigma}$. 
\endproof

\subsection{Convergence of the Lie series operator}
\begin{lem}\label{lem:two}
 Let $d',d'' \in \RR^+$ such that $d'+d'' <1$ and $F,G$ be two functions on $\ml{D}$ such that $\norm{G}{[(1-d')(\rho,\sigma);\zeta]}$ and $\norm{F}{[(1-d'')(\rho,\sigma);\zeta]}$ are bounded for all $\xi \in \ml{R}_{\zeta}$.\\  
 Then, for all $0<d<1-d'-d''$ and all $\nu \in (0,1/2]$, the following inequality holds at each point of $\ml{R}_{\zeta}$ 
\beq{eq:twoparameter}
\norm{\ml{L}_{G} F}{[(1-d-d'-d'')(\rho,\sigma);\zeta]} \leq C \norm{G}{[(1-d')(\rho,\sigma);\zeta]} 
\norm{F}{[(1-d'')(\rho,\sigma);\zeta]}\mx{,}
\eeq
where $C=2[e \rho \sigma (d+d')(d+d'')]^{-1}$.
\end{lem}
\proof
Straightforward\footnote{The different norm used in this paper does not imply substantial differences.} from \cite{gz92}.
\endproof
\begin{prop}\label{prop:chipsi}
Let $d_1,d_2 \in [0,1/2]$ and $\chi$ and $\psi$ be two functions on $\ml{D}$ such that $\norm{\chi}{[(1-d_1)(\rho,\sigma);\zeta]}$ and $\norm{\psi}{[(1-d_2)(\rho,\sigma);\zeta]}$ are bounded for all $\xi \in \ml{R}_{\zeta}$.\\
Then for all $\tilde{d} \in (0,1-\hat{d})$ where $\hat{d}:=\max\{d_1,d_2\}$ and for all $s \geq 1$ one has the following estimate  
\beq{eq:iterativo}
\norm{\ml{L}_{\chi}^s \psi}{[(1-\tilde{d}-\hat{d})(\rho,\sigma);\zeta]} \leq 
\frac{s!}{e^2} \left( \frac{8e}{\rho \sigma \tilde{d}^2}\right)^s
\norm{\chi}{[(1-d_1)(\rho,\sigma);\zeta]}^s 
\norm{\psi}{[(1-d_2)(\rho,\sigma);\zeta]} \mx{.}
\eeq
\end{prop}
\proof
Straightforward going along the lines of Lemma 4.2 of \cite{giorgilli02} and by using\footnote{the factor $8$, in place of $2$ obtained in \cite{giorgilli02}, follows from a rescaling $(\rho,\sigma)\leftarrow (1-\hat{d})(\rho,\sigma)$ and from $\hat{d} \leq 1/2$.} Lemma \ref{lem:two}.
\endproof 

\begin{prop}\label{prop:exp}
In the same hypotheses of Prop. \ref{prop:chipsi}, suppose that, in addition, 
\beq{eq:convergence}
\mathfrak{L}=\frac{8 e}{\tilde{d}^2 \rho \sigma} \norm{\chi}{[(1-d_1)(\rho,\sigma);\zeta]} \leq \frac{1}{2} \mx{.}
\eeq
Then the operator $\exp(\ml{L}_{\chi}) \psi$ is well defined and for all $\tilde{d} \in (0,1-\hat{d})$ the following estimate holds
\beq{eq:estimatelieuno}
\norm{\sum_{s \geq 1} \frac{1}{s!} \ml{L}_{\chi}^s \psi}{[(1-\tild-\hat{d})(\rho,\sigma);\zeta]} \leq \frac{2 \mathfrak{L}}{e^2} \norm{\psi}{[(1-d_2)(\rho,\sigma);\zeta]} \mx{,}
\eeq
in particular
\beq{eq:estimateliedue}
\norm{\exp(\ml{L}_{\chi}) \psi}{[(1-\tild-\hat{d})(\rho,\sigma);\zeta]} \leq 2 \norm{\psi}{[(1-d_2)(\rho,\sigma)]} \mx{.}
\eeq
\end{prop}
 \proof
 It is sufficient to recall the definition of $\exp(\ml{L}_{\chi})$, apply Prop. \ref{prop:chipsi}, and then use $\mathfrak{L} \leq 1/2$.
 \endproof
Note that the previous result holds also if an arbitrary domain restriction $\zeta \rw (1-d)\zeta$ is considered, for all $d \in [0,1)$.
\clearpage
\section{Quantitative estimates on the formal scheme}
Consider the following set of parameters by setting $u_j \equiv (u_j^1,\ldots,u_j^6):=(d_j,\epsilon_j,\zeta_j,m_j,\rho_j,\sigma_j)$ with $u_j^l \in [0,1)$ for all $l=1,\ldots, 6$ and all $j \geq 0$. The vector $u_0$ will be chosen later (see Sec. \ref{seq:inducbasis}).\\  Set, in addition $u_*:=(0,0,0,m_*,\rho_*,\sigma_*)$ for some $m_*,\rho_*,\sigma_*>0$ to be determined (Sec. \ref{sec:controlseq}). As well as for $a$, the property $u_j^l \in [0,1)$ will be repeatedly used in the follow (without an explicit mention) allowing to obtain simpler estimates. 
\begin{lem}\label{lem:iterative}
In the same assumption of Lemma \ref{lem:iterativeformal}, suppose, in addition, the existence of $u_j$ with $u_j>u_*$,  satisfying
\begin{enumerate}
\item 
\beq{eq:iterativeitemone}
\max\left\{\norm{A^{(j)}}{[\sigma_j;\zeta_j]},\norm{B^{(j)}}{[\sigma_j;\zeta_j]}\right\} \leq \epsilon_j e^{-a|\xi|}
\mx{,}
\eeq
\item for all vector valued functions $w=w(q,\xi)$ holds
\beq{eq:iterativeitemthree}
\norm{C^{(j)}(q,\xi) w(q,\xi)}{[\sigma_j;\zeta_j]}  \leq m_j^{-1} \norm{w(q,\xi)}{[\sigma_j;\zeta_j]} \mx{,}
\eeq
\item holds $d_j \leq 1/6$ and $\zeta_j$ is set as
\beq{eq:zetaj}
2 |\omega| \zeta_j=d_j \sigma_j \mx{,}
\eeq
\end{enumerate}
Then there exists a constant $D$ such that: if 
\beq{eq:piccolaunmezzo}
\epsilon_j \frac{D}{a^3 m_j^4 d_j^{4(\tau+1)}} \leq \frac{1}{2} \mx{,}
\eeq
then it is possible to choose $u_{j+1} < u_j$  under the constraint (\ref{eq:zetaj})\footnote{I.e. satisfying $2 |\omega| \zeta_{j+1}=d_{j+1} \sigma_{j+1}$. As well as in the follow, the indices should be changed in $j+1$ where necessary .}, for which (\ref{eq:iterativeitemone}) and (\ref{eq:iterativeitemthree}) are satisfied by $A^{(j+1)},B^{(j+1)}$ and $C^{(j+1)}$ given by (\ref{eq:ajpuno}), (\ref{eq:bjpuno}) and (\ref{eq:cjpuno}), respectively.
\end{lem}
\proof This result is the quantitative counterpart of Lemma \ref{lem:iterativeformal} end this proof is split for the sake of clarity, depending on the considered objects. In order to simplify the notation, the index $j$ will be dropped from all the iterative objects depending on $j$, being restored only in the final estimates.
\subsubsection{Estimates on the generating functions}
Let us consider equation (\ref{eq:firsthomological}). Due to the assumptions, we can apply Prop. \ref{prop:small} with $\delta=0$ and $K=\epsilon$, obtaining
\begin{subequations}
\begin{align}
\norm{\phi}{[(1-d)\sigma; \zeta]} & \leq  
\de \epsilon \frac{M_0}{a d^{\tau} } e^{-a |\xi|} \label{eq:homx} \mx{,} \\
\norm{\pl_q \phi}{[(1-d)\sigma; \zeta]} & \leq \de \epsilon \frac{M_1}{a d^{\tau+1} } e^{-a |\xi|} \label{eq:homxfirst} \mx{,}
\end{align}
\end{subequations}
where $M_0:=S_1\sigma_*^{-\tau}$ and $M_1:=nS_2 \sigma_*^{-(\tau+1)}$.\\
Recalling the definition (\ref{eq:bhat}) then using (\ref{eq:iterativeitemone}), (\ref{eq:homxfirst}) and (\ref{eq:iterativeitemthree}), one gets
\begin{eqnarray}
\norm{\hat{B}}{[(1-d)\sigma;\zeta]} & \leq & 
 \de  \epsilon e^{-a |\xi|} + \frac{1}{m}  \norm{\pl_q \phi}{[(1-d)\sigma;\zeta]} \leq 
\de \epsilon \frac{(1+M_1)}{a m d^{\tau+1}} e^{-a |\xi|} \label{eq:bhatest}\mx{,} \\
\de \norm{\pl_{\xi} \phi}{[(1-d)\sigma;(1-d)\zeta]} & \leq & 
\de \frac{1}{d \zeta} \norm{\phi}{[(1-d)\sigma;\zeta]}
\Heq{\ref{eq:homx}}{\leq} \de  \epsilon \frac{M_0}{a d^{\tau+1} \zeta} e^{-a |\xi|} \label{eq:xxi}\mx{.}
\end{eqnarray}
As for equation (\ref{eq:secondhom}), Prop. \ref{prop:small} used component-wise with $\delta=d$, similarly yields by (\ref{eq:bhatest})
\begin{subequations}
\begin{align}
\norm{Y}{[(1-2d)\sigma;\zeta]} & \leq  
\de \epsilon \frac{M_2 \sigma_*}{a^2 m d^{2 \tau +1}} e^{-a |\xi|} \label{eq:homy} \mx{,} \\
\norm{\pl_q Y}{[(1-2d)\sigma;\zeta]} & \leq \epsilon \frac{M_3}{a^2 m d^{2 \tau +2}} e^{-a |\xi|} \label{eq:homyfirst} \mx{,}
\end{align}
\end{subequations}
where 
\begin{eqnarray}
M_2&:=& \de n S_1 (1+M_1) \sigma_*^{-(\tau+1)} \mx{,}
\label{eq:cdue}\\
M_3&:=& \de n^2 S_2 (1+M_1) \sigma_*^{-(\tau+1)} 
\label{eq:ctre}
\mx{.}
\end{eqnarray}
As a consequence we have, by using (\ref{eq:inequno})
\begin{eqnarray}
\de \norm{\la Y,p \re}{[\rho,(1-2d)\sigma;\zeta]} & \leq & 
\de  \epsilon \frac{M_2 \rho \sigma_*}{a^2 m d^{2 \tau +1}} e^{-a |\xi|} \mx{,} \label{eq:ydotp}\\
\de \norm{Y_{\xi}}{[(1-2d)\sigma;(1-d)\zeta]} & \leq & 
\de \frac{1}{d \zeta} \norm{Y}{[(1-2d)\sigma;\zeta]} 
 \leq 
\de \epsilon \frac{M_2}{a^2 m d^{2 \tau +2} \zeta} e^{-a |\xi|} \label{eq:yxi}\mx{.}
\end{eqnarray}
By (\ref{eq:ydotp}), Prop. \ref{prop:exp} and setting $\mathfrak{L}:=Q_1 e^{-a|\xi|}$,  we have that $\exp(\ml{L}_{\la Y,p \re})$ converges uniformly on $\ml{R}_{\zeta}$ provided\footnote{In this case $d_1:=2d$, while $d_2 \leq 2d$ as used below, so it is possible to set $\tilde{d} \equiv d <1-2d$ by hypothesis $(3)$. Moreover, the latter implies $d_1,d_2 \leq 1/2$ as required by Prop. \ref{prop:chipsi}.}
\beq{eq:convergenceliesectre}
Q_1:=\epsilon \frac{8 e M_2}{a^2 m d^{2 \tau+3}} \leq \frac{1}{2}
\eeq
\subsubsection{Estimates on the transformed Hamiltonian}
Firstly, by (\ref{eq:ahat}), using (\ref{eq:iterativeitemthree}) and (\ref{eq:homxfirst}) one gets 
\[ \norm{\hat{A}}{[(1-d)\sigma;\zeta]} 
\leq \epsilon^2 \frac{M_1(1+M_1)}{a^2 m d^{2 \tau+2}} e^{-2a |\xi|}\mx{.}
\] 
 Hence by (\ref{eq:ajpuno}), Prop. \ref{prop:exp} with $d_2=d$ and after an arbitrary restriction in $\rho$ and $\zeta$, we have
 \beq{eq:stimaaprimo}
 \norm{A^{(j+1)}}{[(1-3d_j)(\rho_j,\sigma_j;\zeta_j)]}  \leq   \de \epsilon_j^2 \frac{M_4}{a^2 m_j d_j^{2 \tau+2}} e^{-2 a|\xi|}  \mx{,}
 \eeq
  where
 \beq{eq:cquattro}
 M_4:= 2 M_1(1+M_1) \mx{.}
 \eeq
 On the other hand, by (\ref{eq:estimatelieuno}),  (\ref{eq:bhatest}) and (\ref{eq:inequno}) 
\[
 \de \norm{\sum_{s \geq 1} \frac{s}{(s+1)!} \ml{L}_{\la Y,p \re}^s \la \hat{B}, p \re}{[(1-3d)(\rho,\sigma);\zeta]}  \leq   \de \frac{2 \mathfrak{L}}{e^2} \norm{\la \hat{B},p \re}{[(1-d)(\rho,\sigma);\zeta]} \leq
   \epsilon \frac{2 \rho (1+M_1)Q_1}{a m e^2 d^{\tau+1}} e^{-2a|\xi|} \mx{.}
\]
Recalling (\ref{eq:bjpuno}), the definition in (\ref{eq:convergenceliesectre}) and (\ref{eq:inequnocontrary}), 
\beq{eq:stimabprimo}
 \de \norm{ B^{(j+1)}}{[(1-3d_j)(\rho_j,\sigma_j;\zeta_j)]} \leq \epsilon_j^2 \frac{M_5}{a^3 m_j^2 d_j^{3 \tau+4}}e^{-2 a |\xi|} \mx{,}
\eeq
with 
\beq{eq:ccinque}
M_5:=16 n (1+M_1) M_2(e \sigma_*)^{-1} \mx{.}
\eeq
Let us set $C':=C^{(j+1)}$. Directly from (\ref{eq:cjpuno}), Prop. \ref{prop:exp} and (\ref{eq:iterativeitemthree}) one has 
\beq{eq:stimaccprimo}
\norm{\la (C'-C)p,p \re}{[(1-3d)(\rho,\sigma);\zeta]} \leq \frac{2 \mathfrak{L}}{e^2} \norm{\la Cp,p \re}{[(1-2d)(\rho,\sigma);\zeta]} \leq \epsilon \frac{16 M_2 }{a m^3 e d^{2 \tau+3}} \rho^2 e^{-a|\xi|} \mx{,}
\eeq
implying, by (\ref{eq:ineqdue})
\beq{eq:estimateckl}
\norm{C_{kl}'-C_{kl}}{[(1-3d)\sigma;\zeta]} \leq \epsilon \frac{M_6 }{a^2 m^3 n d^{2 \tau+3}} e^{-a|\xi|} \mx{,}
\eeq
with
\beq{eq:csei}
M_6:=16 n M_2(e \sigma_*)^{-1} \mx{.}
\eeq
Now set
\beq{eq:mprimo}
m':=m- \epsilon \frac{M_6 }{a^2 m^3 d^{2 \tau+3}} e^{-a|\xi|}\mx{,}
\eeq
which is well defined provided that, e.g.
\beq{eq:limcsei}
\epsilon \frac{M_6 }{a^2 m^4 d^{2 \tau+3}} \leq \frac{1}{2} \mx{.}
\eeq
giving, in particular, $m' \in [m/2,m]$. In this way we have for all $w=w(q,\xi)$
\beq{eq:stimacprimo}
\begin{array}{rcl}
\de \norm{C'w}{[(1-3d)\sigma;\zeta]} & \Heq{\ref{eq:iterativeitemthree})(\ref{eq:stimaccprimo}}{\leq} &  
\de \left(\frac{1}{m}+ \epsilon \frac{M_6 }{a^2 m^3 d^{2 \tau+3}} e^{-a|\xi|} \right) \norm{w}{[(1-3d)\sigma;\zeta]}  \\
& \leq & \de \frac{1}{m'} \norm{w}{[(1-3d)\sigma;\zeta]} \mx{,}
\end{array}
\eeq
where the inequality $a^{-1}+b<(a-b)^{-1}$, valid for all $0<b<a<1$, then (\ref{eq:mprimo}) have been used in the last passage.
\subsubsection*{Determination of parameters} Let us set
\beq{eq:epjpuno}
\epsilon_{j+1}:=\frac{D}{a^3 m_j^4 d_j^{4(\tau+1)}} \epsilon_j^2 \mx{.}
\eeq
In this way, conditions (\ref{eq:convergenceliesectre}), (\ref{eq:limcsei}) and those obtained by comparing (\ref{eq:stimaaprimo}) and (\ref{eq:stimabprimo}) with (\ref{eq:iterativeitemone}), are implied \emph{a fortiori} by hypothesis (\ref{eq:piccolaunmezzo}), provided that $
D:=\max\{8 e M_2 ,M_4 , M_5, M_6\}$. The property $\epsilon_{j+1}<\epsilon_j$ is an easy consequence of (\ref{eq:piccolaunmezzo}) and of $\epsilon_j<1$.\\
By taking into account the estimates (\ref{eq:stimaaprimo}) and (\ref{eq:stimabprimo}), we have that the domain on which these hold requires the restriction described by the following choices 
\beq{eq:sigmarhojpuno}
\sigma_{j+1}:=(1-3 d_j) \sigma_j, \qquad \rho_{j+1}:=(1-3 d_j) \rho_j \mx{.}
\eeq
As for $\zeta_{j+1}$, condition (\ref{eq:sceltazeta}) is valid at the $j+1-$th step if 
$
\zeta_{j+1}=(2 |\omega|)^{-1} \min \{ (1-3 d_j)d_j \sigma_j, d_{j+1} \sigma_{j+1}\}
$. 
As $d_j \leq 1/6$ by hypothesis, by the first of (\ref{eq:sigmarhojpuno}) the previous condition is of the form 
(\ref{eq:zetaj}) provided that $d_{j+1} < d_j$ is chosen. This implies $\zeta_{j+1}<\zeta_j$.\\
The only parameter left is $m_j$. Note that  (\ref{eq:piccolaunmezzo}) implies, in particular
$\epsilon M_6/(a^2 m^3 d^{2 \tau+3})  \leq m d^{2 \tau+1} $, then
\[
m':=m-\epsilon \frac{M_6}{a^2 m^3 d^{2 \tau+3} \zeta} e^{-a|\xi|} \geq m(1-d^{2\tau+1}) \mx{.}
\]
In conclusion, inequality (\ref{eq:stimacprimo}), hence (\ref{eq:iterativeitemthree}), are satisfied by setting 
\beq{eq:mjpuno}
m_{j+1}:=m_j(1-d_j^{2\tau+1})  \mx{.}
\eeq
The choice of $u_{j+1}$ is now complete\footnote{The freedom in the choice of $d_{j+1}$ (subject only to the constraint $d_{j+1}<d_j$) will be profitably used later.}.
\endproof

\subsection{Estimates on the transformation of variables}

\begin{prop}\label{prop:trasf} Assume the validity of Lemma \ref{lem:iterative}. Then, for all $j \in \NN$, the transformation (\ref{eq:change}) is a symplectic transformation
\[
\ml{K}_j: \ml{D}_{j+1} \longrightarrow \ml{D}_j \mx{,}
\]
where $\ml{D}_j:=\Delta_{\rho_j}(0)\times \TT_{2 \sigma_j}^n \x \ml{S}_{\rho_j} \x \ml{R}_{\zeta_j} \ni (q^{(j)},p^{(j)},\eta^{(j)},\xi^{(j)})$, for which there exists a constant $T$ such that, 
\begin{subequations}
\begin{align}
|q^{(j+1)}-q^{(j)}| & \leq T \sigma_j d_j e^{-a|\xi|} \mx{,}\\
|p^{(j+1)}-p^{(j)}| & \leq T \rho_j d_j e^{-a|\xi|} \mx{,}\\
|q^{(j+1)}-q^{(j)}| & \leq T \rho_j d_j e^{-a|\xi|} \mx{,}
\end{align}
while $|\xi^{(j+1)}-\xi^{(j)}|=0$, i.e. $\xi^{(j)}=:\xi$ for all $j$. Moreover $\ml{K}_j$ is $\epsilon_0-$``close to the identity'', i.e. $\lim_{\epsilon_0 \rw 0} \ml{K}_j=\id$ for all $j$.
\end{subequations}
\end{prop}
\proof
Once more it is convenient to examine separately the transformations realising $\ml{K}_j$
\begin{align*}
(\hat{q}^{(j)},\hat{p}^{(j)},\hat{\eta}^{(j)},\hat{\xi}^{(j)})&:=
\exp(\ml{L}_{\phi^{(j)}}) (q^{(j+1)},p^{(j+1)},\eta^{(j+1)},\xi^{(j+1)}) \mx{,}\\ 
(q^{(j)},p^{(j)},\eta^{(j)},\xi^{(j)})&:=
\exp(\ml{L}_{\chi^{(j)}})
(\hat{q}^{(j)},\hat{p}^{(j)},\hat{\eta}^{(j)},\hat{\xi}^{(j)}) \mx{.}
\end{align*}
Due to the structure of $\phi^{(j)}$ the action of the first operator reduces to the first term for the momenta, 
\begin{align*}
\hat{p}^{(j)}&=p^{(j+1)}+ [\pl_{q} \phi^{(j)}]_{(q,\xi)=(q^{(j+1)},\xi^{(j+1)})} \mx{,}\\ 
\hat{\eta}^{(j)}&=\eta^{(j+1)}+[\pl_{\xi} \phi^{(j)}]_{(q,\xi)=(q^{(j+1)},\xi^{(j+1)})} \mx{,}
\end{align*}
while it is the identity in the other variables: $\hat{q}^{(j)}=q^{(j+1)}$ and $\hat{\xi}^{(j)}=\xi^{(j+1)}$. Quantitatively we find
\[
|\hat{p}^{(j)}-p^{(j+1)}| \Heq{\ref{eq:homxfirst}}{\leq} 
\epsilon_j \frac{M_1}{a d_j^{\tau+1}} e^{-a |\xi^{(j+1)}|},\qquad 
|\hat{\eta}^{(j)}-\eta^{(j+1)}| \Heq{\ref{eq:xxi}}{\leq} 
\epsilon_j 
\frac{M_0}{a d_j^{\tau+1} \zeta_j} e^{-a |\xi^{(j+1)}|} \mx{.}
\]
As for the second transformation, first note that 
\beq{eq:four}
\ml{L}_{\chi^{(j)}} q =  Y^{(j)}, \quad 
\ml{L}_{\chi^{(j)}} p  =  \la \pl_{q} Y^{(j)} ,p \re  \quad
\ml{L}_{\chi^{(j)}} \xi  =   0 ,\quad
\ml{L}_{\chi^{(j)}} \eta  =  \la  \pl_{\xi} Y^{(j)}, p \re \mx{,}
\eeq
where the expressions above are meant to be evaluated at $(q,p,\eta,\xi)=
(\hat{q}^{(j)},\hat{p}^{(j)},\hat{\eta}^{(j)},\hat{\xi}^{(j)})$. Now consider bound (\ref{eq:iterativo}) for $s-1$, setting $\chi:=\chi^{(j)}$ and $\psi$ as the objects in the (\ref{eq:four}) r.h.sides one by one. We get, e.g., for the first of them 
\[
\norm{\ml{L}_{\chi^{(j)}}^s q}{[(1-3d_j)(\rho_j,\sigma_j;\zeta_j)]}  \leq   
\frac{s!}{e^2} \mathfrak{L}^{s-1} 
\norm{Y^{(j)}}{[\rho_j,(1-2d_j)\sigma_j;(1-d_j)\zeta_j]} \leq 
s! \frac{d^2 \sigma_* }{8 e^3} \mathfrak{L}^s \mx{.}
\]
Repeating this computation also for the other variables we get (recall $\sum_{s \geq 1} \mathfrak{L}^s \leq 2 \mathfrak{L}$) 
\begin{subequations}
\begin{align}
\de |q^{(j+1)}-\hat{q}^{(j)}| & \leq  \de \frac{d_j^2 \sigma_*}{4 e^3} \mathfrak{L} = \epsilon_j \frac{2 M_2}{a^2 e^2 m_j  d_j^{2 \tau+1}} e^{-a |\xi^{(j)}|} \label{eq:qhatal} \mx{,}\\
\de |p^{(j+1)}-\hat{p}^{(j)}| & \leq  \de \frac{ d_j \rho_j M_3}{4 e^3 M_2} \mathfrak{L}
= \epsilon_j \frac{2 M_3 \rho_j}{a^2 e^2 m_j  d_j^{2 \tau+2} } e^{-a |\xi^{(j)}|} \mx{,}\\
\de |\eta^{(j+1)}-\hat{\eta}^{(j)}| & \leq  \de \frac{d_j \rho_j}{4 e^3 \zeta_j} \mathfrak{L}
= \epsilon_j \frac{2 M_2 \rho_j}{a^2 e^2 m_j  d_j^{2 \tau+2} \zeta_j} e^{-a |\xi^{(j)}|} \mx{,}
\end{align}
\end{subequations}
and clearly $\de |\xi^{(j+1)}-\hat{\xi}^{(j)}|  =   0$, implying $\xi^{(j+1)} \equiv \xi^{(j)}$.
\begin{rem} It is finally evident that the transformation $\ml{K}_j$ does not act on time, hence we can set $\xi^{(j)} \equiv \xi$ for all $j \in \NN$ as in the statement. On the other hand this is a necessary property in order to obtain a meaningful result.
\end{rem}
Collecting the obtained estimates we get that $|q^{(j+1)}-q^{(j)}|$ is given by (\ref{eq:qhatal}), while
\beq{eq:pieta}
\begin{array}{rcl}
\de |p^{(j+1)}-p^{(j)}| & \leq &
\de \epsilon_j  \frac{(M_1 e^2+2 M_3) \rho_j}{a^2 e^2 m_j d_j^{2\tau+2}  \rho_*} e^{-a|\xi|} \mx{,}
\\
\de |\eta^{(j+1)}-\eta^{(j)}| & \leq & \epsilon_j \de   \frac{(M_0 e^2+2 M_2) \rho_j}{a^2 e^2 m_j d_j^{2\tau+2} \rho_*  \zeta_j} e^{-a|\xi|} \mx{,}
\end{array}
\eeq
having used $\rho_j>\rho_*$. Hence it is possible to find\footnote{Precisely $T:=(D e^2 \rho_* \sigma_*)^{-1}\max\{ M_2 \rho_*,(M_1 e^2+2 M_3)\sigma_*,2|\omega|(M_0 e^2+2M_2)\}$, by (\ref{eq:qhatal}), (\ref{eq:pieta}) and using (\ref{eq:piccolaunmezzo}) and (\ref{eq:zetaj}).} $T$, obtaining the desired estimates.\\
The $\epsilon_0-$closeness to the identity easily follows from (\ref{eq:qhatal}), (\ref{eq:pieta}) and from the monotonicity of $\{\epsilon_j\}$.
\endproof
\clearpage
\section{Convergence of the formal scheme}
\subsection{Construction of the control sequence}\label{sec:controlseq}
\begin{lem} In the assumptions of Lemma \ref{lem:iterative}, it is possible to determine $u_*$ and construct the sequence $\{u_j\}_{j \in \NN}$ such that 
\beq{eq:limit}
\lim_{j \rw \infty} u_j=u_* \mx{.}
\eeq
\end{lem}
\proof
Let us choose in (\ref{eq:epjpuno}) $\epsilon_j=\epsilon_0 j^{-8(\tau+1)}$, obtaining
\beq{eq:dj}
d_j  =  \left(\frac{D \epsilon_0}{a^3 m_j^4} \right)^{\frac{1}{4(\tau+1)}} \frac{(j+1)^2}{j^4} \mx{.}
\eeq
The following bound is immediate for all $j \geq 1$
\beq{eq:bounddk}
d_j \leq \de 2  \frac{\ml{A}}{j^2},\qquad \ml{A}:=\left(\frac{D \epsilon_0}{a^3 m_*^4} \right)^{\frac{1}{4(\tau+1)}} \mx{.}
\eeq
Imposing condition $d_{j} \geq d_{j+1}$ in (\ref{eq:dj}) one gets $(1-d_j^{2 \tau+3})^{\frac{1}{\tau+1}} \geq j^4(j+2)^2/(j+1)^6$. By using (\ref{eq:bounddk}), it takes the stronger form
\[
1-2 \ml{A} j^{-2} \geq \frac{j^4(j+2)^2}{(j+1)^6} \mx{.}
\]
The latter is true for all $j$ provided that it holds for $j=1$. This is achieved if $\ml{A} \leq 55/128$, a condition that can be enforced by requiring  $\ml{A} \leq 1/12$. In this way we obtain $d_j \leq d_1  \leq 1/6$ as required by Lemma \ref{lem:iterative}, item ($3$). This immediately implies
\beq{eq:seriesdj}
\sum_{j \geq 1} d_j \leq \frac{1}{6} \sum_{j \geq 1} j^{-2} < \left(\frac{\pi}{6}\right)^2 \mx{.} 
\eeq
In this way, the range of the admissible values for $\epsilon_0$ is determined once and for all; more explicitly
\beq{eq:limep}
\frac{D \epsilon_0}{a^3 m_*^4}  \leq \frac{1}{12^{4(\tau+1)}} \mx{.}
\eeq
We only need to prove the limit (\ref{eq:limit}). Let us start from $\rho_j$. By (\ref{eq:sigmarhojpuno}) we have that if $\prod_{j \geq 1} (1-3 d_j)$ 
is lower bounded by a constant, say $M_{\rho}$, then $\rho_0 M_{\rho} $ is a lower bound for $\rho_j$ for all $j$.\\
Consider
\[
\log \prod_{j \geq 1} (1-3 d_j) = \sum_{j \geq 1} \log (1-3 d_j) \geq - 6  \log 2 \sum_{j \geq 1} d_j >  -\log 4 \mx{,}
\]
in which we have used the inequality $0 \geq \log(1-x) \geq -2 x \log 2$, valid for $x \in [0,1/2]$. Hence $\prod_{j \geq 1} (1-3 d_j) \leq 1/4$. This implies that the required lower bound holds for $\rho_*=\rho_0/4$ and then $\sigma_*=\sigma_0/4$.
A similar arguments applies for $m_j$, yielding  $m_*=m_0/2$.
\endproof

\subsection{Induction basis and conclusion of the proof}\label{seq:inducbasis}
In this final part we check that the inductive hypotheses described in Lemmas \ref{lem:iterativeformal} and \ref{lem:iterative} hold at the initial step, i.e. $j=0$, fixing in this way $u_0$.\\
First of all we see that $H$ is of the form (\ref{eq:hamricorsiva}) in a way we can set $H_0:=H$. It is sufficient to consider the (finite) Taylor  expansion of $f$ around $p=0$ in (\ref{eq:ham}) then define
\[
A^{(0)}:=\ep f(q,0,\xi),\qquad 
B^{(0)}:=\ep \pl_p f(q,0,\xi),\qquad 
C^{(0)}:=\Gamma+\ep \pl_p^2 f(q,0,\xi) \mx{.}
\]
Note that $C^{(0)}$ is symmetric. Now set $\rho_0:=\rho/2$ and $\sigma_0:=\sigma$. By a Cauchy estimate and (\ref{eq:slowdecay}) we have 
\beq{eq:lastcauchy}
\norm{\pl_p f}{[\rho_0,\sigma_0;\zeta_0]} \leq 
M_f \rho_0^{-1} e^{-a |\xi|},\qquad \norm{\pl_p^2 f}{[\rho_0,\sigma_0;\zeta_0]} \leq 
 M_f \rho_0^{-2}e^{-a |\xi|} \mx{,}
\eeq
for all $\zeta_0$ (determined below). Hence (\ref{eq:iterativeitemone}) is satisfied for $j=0$ by setting $\epsilon_0:=\ep M_f / \rho_0$. By Prop. \ref{prop:trasf}, this shows that the sequence $\{\ml{K}_j\}$ and then the composition
\beq{eq:composition}
\ml{K}:=\lim_{j \rw \infty} \ml{K}_j \circ \ml{K}_{j-1} \circ \ldots \circ \ml{K}_0 \mx{,}
\eeq
is $\ep-$close to the identity.\\It is natural to realize that (\ref{eq:iterativeitemthree}) holds by virtue of (\ref{eq:hypongamma}) and for sufficiently small $\ep$. From the quantitative point of view one can ask $ |C^{(0)} v| \leq m_0^{-1}|v|$ for all $v \in \CC^n$ with $m_0:=m/2$. This is true for all $\ep \leq \tilde{\ep}$ where
\beq{eq:limitepzero}
\tilde{\ep}:=\rho^2 (16 M_f n)^{-1}(\sqrt{m^2 \norm{\Gamma}{\infty}^2+12}-m \norm{\Gamma}{\infty}) \mx{,}
\eeq
denoted\footnote{This bound follows from a straightforward check. By the second of (\ref{eq:lastcauchy}) we have $C^{(0)}=\Gamma+\ep h E$ where $h:= M_f \rho_0^{-2}$ and $E_{kl} \in [-1,1]$ for all $k,l=1,\ldots,n$. It is now sufficient to use the (exact) Mac Laurin expansion $|C^{(0)}(\ep)v|^2=|\Gamma v|^2+2 \ep h \la \Gamma v, E v \re + \ep^2 |E v|^2$ and (\ref{eq:hypongamma}) to get the (\ref{eq:limitepzero}).}  
$\norm{\Gamma}{\infty}:=\max_i \sum_{j=1}^n |\Gamma_{ij}|$.\\ 
The choice of $u_0$ is now complete by choosing $d_0=1/6$ and $\zeta_0$ as determined by (\ref{eq:sceltazeta}). By using (\ref{eq:limep}) and recalling the choice for $\epsilon_0$ and $m^*$ above, we finally obtain the limitation for $\ep_a$
\beq{eq:finalvalueep}
\ep_a = \min\{ \rho a^3 m^4 (2^9 12^{4(\tau+1)} D M_f)^{-1},\tilde{\ep}\} \mx{.}
\eeq
The validity\footnote{The allowed range for $\ep$ found above, exploits a well known issue in the KAM theory: the numerical coefficient in 
(\ref{eq:finalvalueep}) is smaller than $10^{-9}$ (and rapidly decreasing as the number $n$, hence $\tau$, increase). This value is practically unsuitable for interesting physical applications (such as Celestial Mechanics problems).  A relevant branch of the KAM theory is devoted to the development of tools capable to increase this threshold. See \cite{celgiorloc} for an example or \cite{cellettikam} for a comprehensive application of the computer-assisted proofs approach.} of condition (\ref{eq:piccolaunmezzo}) for $j=0$ follows from (\ref{eq:limep}).\\ 
The very last step consists in showing the convergence of the composition (\ref{eq:composition}). By Prop. \ref{prop:trasf} and recalling (\ref{eq:seriesdj}) we find
\[
|q_{\infty}-q| \leq T \sum_{k \geq 0} |q_{k+1}-q_k| < 2 \sigma T \mx{.}
\]
Analogously we find $|p_{\infty}-p|,|\eta_{\infty}-\eta| < 2 \rho T $. Hence by the Weierstra{\ss} Theorem (see, e.g. \cite{dettman}) the transformation (\ref{eq:composition}) converges uniformly in all compact subsets of $\ml{E}_*:=\Delta_{\rho_*} \times \TT_{2 \sigma_*}^n \times \ml{S}_{\rho_*}$. Note that the degeneration of $\ml{R}_{\zeta_j}$ is not an issue as the transformation is trivial in the $\xi$ variable. The proof is completed by setting $\ml{D}_*= \ml{E}_* \times \RR^+$.

\subsection*{Acknowledgements} We would like to acknowledge useful e-mail exchanges with Prof. Antonio Giorgilli.

\bibliographystyle{alpha}
\bibliography{KAM.bib}

\end{document}